\documentclass[12pt]{article}
\usepackage{amsmath,amssymb,amsthm}
\usepackage{txfonts,bm,pifont}
\usepackage{graphicx}
\usepackage[pdfencoding=auto,bookmarks=true,bookmarksnumbered=true,colorlinks=true]{hyperref}
\usepackage{xcolor}
\usepackage{float}
\usepackage{multirow,bigdelim}
\theoremstyle{definition}
\newtheorem{dfn}{Definition}
\theoremstyle{plain}
\newtheorem{thm}{Theorem}
\theoremstyle{remark}
\newtheorem{rmk}{Remark}
\newcommand{\N}{\mathbb{N}}
\newcommand{\Z}{\mathbb{Z}}

\newcommand{\C}{\mathbb{C}}
\newcommand{\Log}{\mathrm{Log}}
\graphicspath{{figs/}}

\usepackage{fancyhdr}
\usepackage{caption}
\begin{document}
\begin{center}
\begin{large}
    Michell--Prager type truss structures constructed from integrable discrete power function and discrete logarithmic function\\[5mm]
\end{large}
\begin{normalsize}
    Motoki \textsc{Masada}\\
    Joint Graduate School of Mathematics for Innovation, Kyushu University\\
    744 Motooka, Nishi-ku, Fukuoka 819-0382, Japan\\
    e-mail: masada.motoki.382@s.kyushu-u.ac.jp\\[2mm]
    Kentaro \textsc{Hayakawa}\\
    College of Industrial Technology, Nihon University\\
    1-2-1 Izumi-cho, Narashino, Chiba 275-8575, Japan\\
    e-mail: hayakawa.kentaro@nihon-u.ac.jp\\[2mm]
    Kazuki \textsc{Hayashi}\\
    Graduate School of Engineering, Kyoto University\\
    Kyoto daigaku-katsura, Nishikyo-ku, Kyoto, Japan 615-8540\\
    e-mail: hayashi.kazuki@archi.kyoto-u.ac.jp\\[2mm]
    Yoshiki \textsc{Jikumaru}\\
    Faculty of Information Networking for Innovation and Design, Toyo University\\
    1-7-11 Akabanedai, Kita-ku, Tokyo, 115-8650, Japan\\
    e-mail: jikumaru@toyo.jp\\[2mm]
    Kenji \textsc{Kajiwara}\\
    Institute of Mathematics for Industry, Kyushu University\\
    744 Motooka, Fukuoka 819-0395, Japan\\
    e-mail: kaji@imi.kyushu-u.ac.jp\\[2mm]
    Yohei \textsc{Yokosuka}\\
    Graduate School of Science and Engineering, Kagoshima University\\
    1-21-40 Korimoto Kagoshima-city Kagoshima, 890-0065, Japan\\
    e-mail: yokosuka@aae.kagoshima-u.ac.jp
\end{normalsize}
\end{center}
\begin{abstract}
    The Michell--Prager type truss structures are constructed from the integrable discrete power and logarithmic functions.
    It is demonstrated that specific sublattices, subject to suitable boundary conditions, yield approximate Michell trusses, which are theoretically optimal structures achieving force equilibrium with material economy.
    The result of shape optimization minimizing the Michell functional is provided for numerical evidence of their optimality.
\end{abstract}


\section{Introduction}
In structural design, finding an ``optimal'' configuration is a fundamental challenge.
The Michell truss is a recognized theoretically ideal structure achieving force equilibrium with material economy.
However, since the original theory assumes members form a continuous orthogonal curvilinear network, constructing practical discrete truss approximations remains a non-trivial task.

The discrete isothermic net, introduced by Bobenko and Pinkall\cite{bobenko1996isothermic} as a discretization of the important class of isothermic surfaces, is expected to be compatible with Michell theory since the orthogonality of member arrangements is a key requirement for optimality.
Recently, it has been demonstrated in\cite{hayashi2024parametric} that truss structures in equilibrium can be constructed from discrete isothermic nets, showing that the truss based on the discrete exponential function effectively approximates the Michell truss.

In this paper, we extend this approach to investigate truss structures constructed from the discrete power and logarithmic functions.
Shape optimization numerically demonstrates that these structures are promising Michell truss approximations.
The results suggest that the geometric properties of discrete isothermic nets are inherently suited for approximating optimal forms.

\section{Discrete holomorphic function}
Bobenko and Pinkall introduced a definition for the integrable discretization of holomorphic functions.
For a map $f : \Z^2 \to \C$ and a fixed $(m, n) \in \Z^2$, we use the following notation: $f=f(m, n)$, $f_{\pm 1}=f(m \pm 1, n)$, $f_{\pm 2}=f(m, n \pm 1)$, $f_{\pm 1, \pm 2}=f(m \pm 1, n \pm 1)$.
\begin{dfn}[Bobenko--Pinkall\cite{bobenko1996isothermic}]
    A map $f : \Z^2 \to \C$ is called a \emph{discrete holomorphic function} if $f$ satisfies
    \begin{equation}
        \frac{(f-f_1)(f_{12}-f_2)}{(f_1-f_{12})(f_2-f)}=-1.
        \label{eq: cr-eq}
    \end{equation}
    \label{def: d_hol}
\end{dfn}
\begin{rmk}
    A map $f : \Z^2 \to \C$ is called a \emph{circular net} if the right hand side of \eqref{eq: cr-eq} is real, which means the four points are concircular.
\end{rmk}

\begin{dfn}[Bobenko\cite{bobenko1999discrete}]
    Let $0<\gamma<2$ and $\Z^2_+:=\{(m, n) \in \Z^2 \mid m \geq 0, n \geq 0\}$.
    A map $Z^\gamma : \Z^2_+ \to \C$ is called the \emph{discrete power function} with the exponent $\gamma$ if $f=Z^\gamma$ satisfies \eqref{eq: cr-eq} and
    \begin{equation}
        \gamma f=2m\frac{(f_1-f)(f-f_{-1})}{f_1-f_{-1}}+2n\frac{(f_2-f)(f-f_{-2})}{f_2-f_{-2}},
        \label{eq: sim_pow}
    \end{equation}
    with initial conditions
    \begin{equation}
        f(0, 0)=0, \quad f(1, 0)=1, \quad f(0, 1)=e^{\gamma\pi i/2}.
    \end{equation}
\end{dfn}
\begin{dfn}[Agafonov--Bobenko\cite{agafonov2000discrete}]
    \emph{Discrete power function} $Z^2 : \Z^2_+ \to \C$ is the solution of \eqref{eq: cr-eq} and \eqref{eq: sim_pow} with initial conditions
    \begin{equation}
        \begin{gathered}
            f(0, 0)=f(1, 0)=f(0, 1)=0,\\
            f(2, 0)=1, \quad f(0, 2)=-1, \quad f(1, 1)=\frac{2}{\pi}i.
        \end{gathered}
    \end{equation}
\end{dfn}
\begin{dfn}[Agafonov--Bobenko\cite{agafonov2000discrete}]
    \emph{Discrete logarithmic function} $\Log : \Z^2_+ \to \overline{\C}$ is the solution of \eqref{eq: cr-eq} and
    \begin{equation}
        1=m\frac{(f_1-f)(f-f_{-1})}{f_1-f_{-1}}+n\frac{(f_2-f)(f-f_{-2})}{f_2-f_{-2}},
        \label{eq: sim_log}
    \end{equation}
    with initial conditions
    \begin{equation}
        \begin{gathered}
            f(0, 0)=\infty, \quad f(1, 0)=0, \quad f(0, 1)=\pi i,\\
            f(2, 0)=1, \quad f(0, 2)=1+\pi i, \quad f(1, 1)=\frac{\pi}{2}i.
        \end{gathered}
    \end{equation}
\end{dfn}
\begin{figure}[htbp]
    \centering
    \begin{minipage}[b]{0.3\linewidth}
        \centering
        \includegraphics[height=4cm]{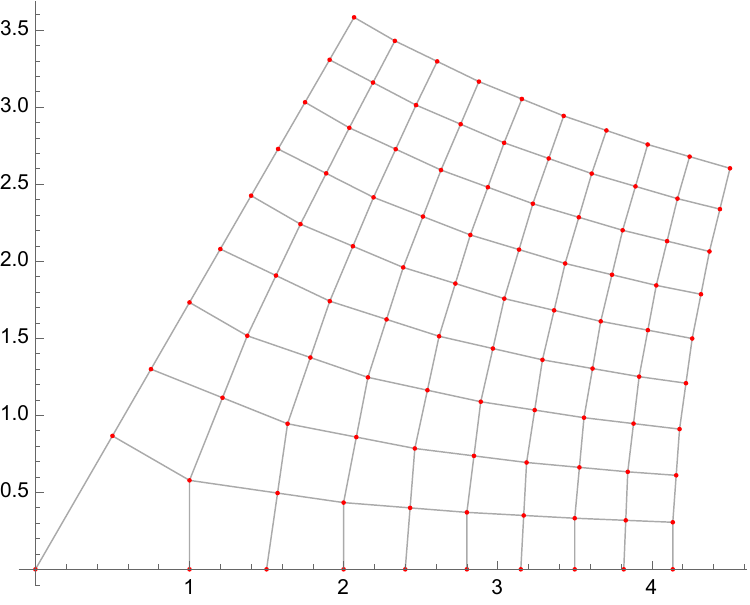}
    \end{minipage}
    \hfill
    \begin{minipage}[b]{0.3\linewidth}
        \centering
        \includegraphics[height=4cm]{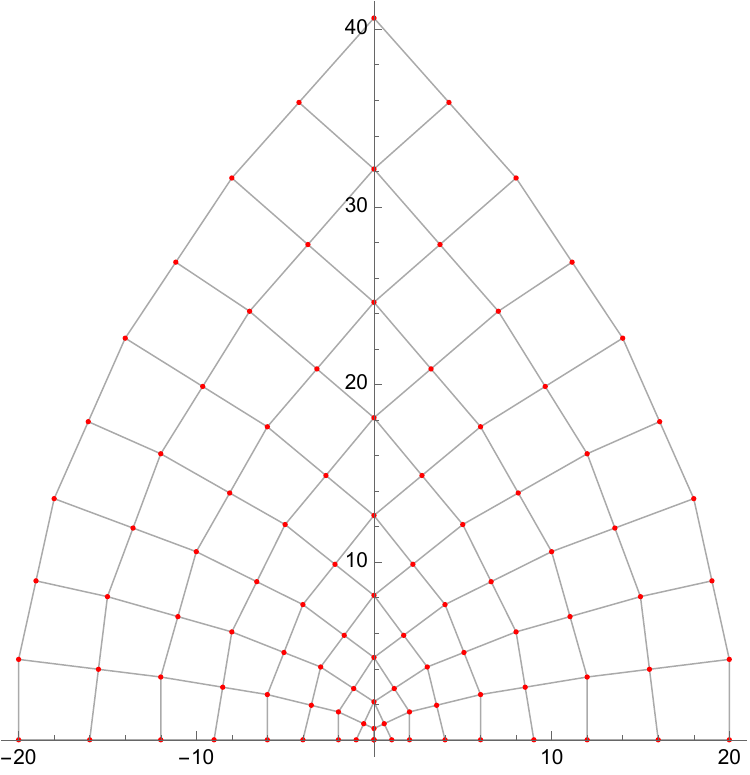}
    \end{minipage}
    \hfill
    \begin{minipage}[b]{0.3\linewidth}
        \centering
        \includegraphics[height=3.5cm]{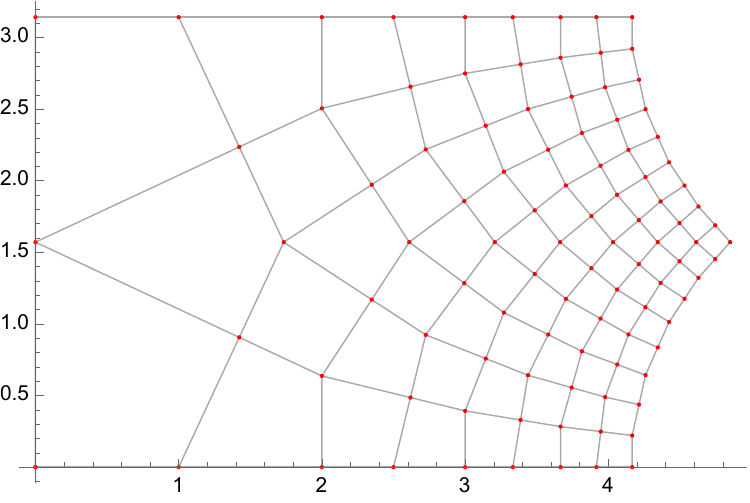}
    \end{minipage}
    \caption{The discrete power functions $Z^\frac{2}{3}$ (left) and $Z^2$ (center), and the discrete logarithmic function $\Log$ (right).}
\end{figure}

\begin{figure}[htbp]
    \centering
    \begin{minipage}[b]{0.45\linewidth}
        \centering
        \includegraphics[height=4.5cm]{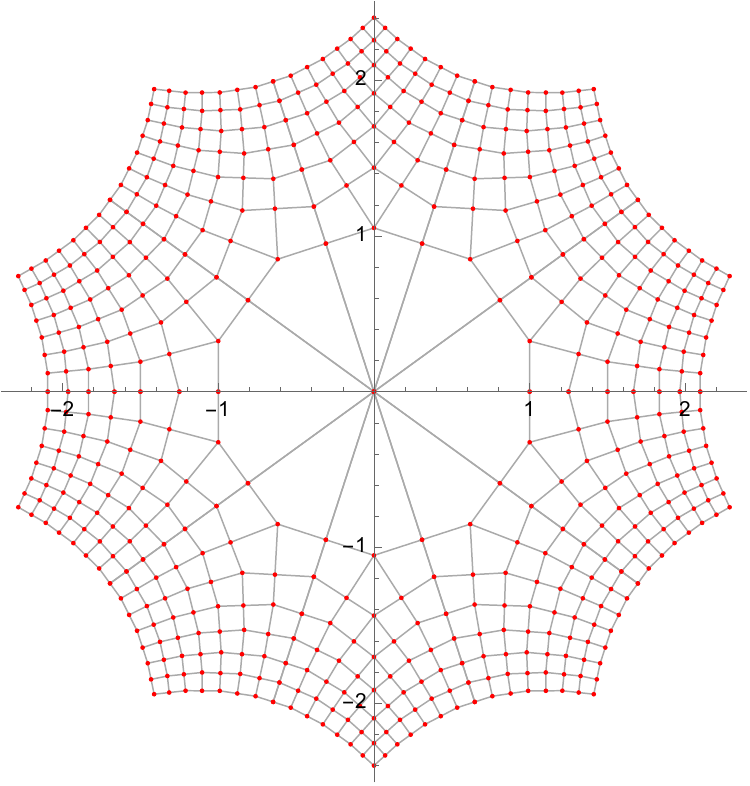}
    \end{minipage}
    \begin{minipage}[b]{0.45\linewidth}
        \centering
        \includegraphics[height=4.5cm]{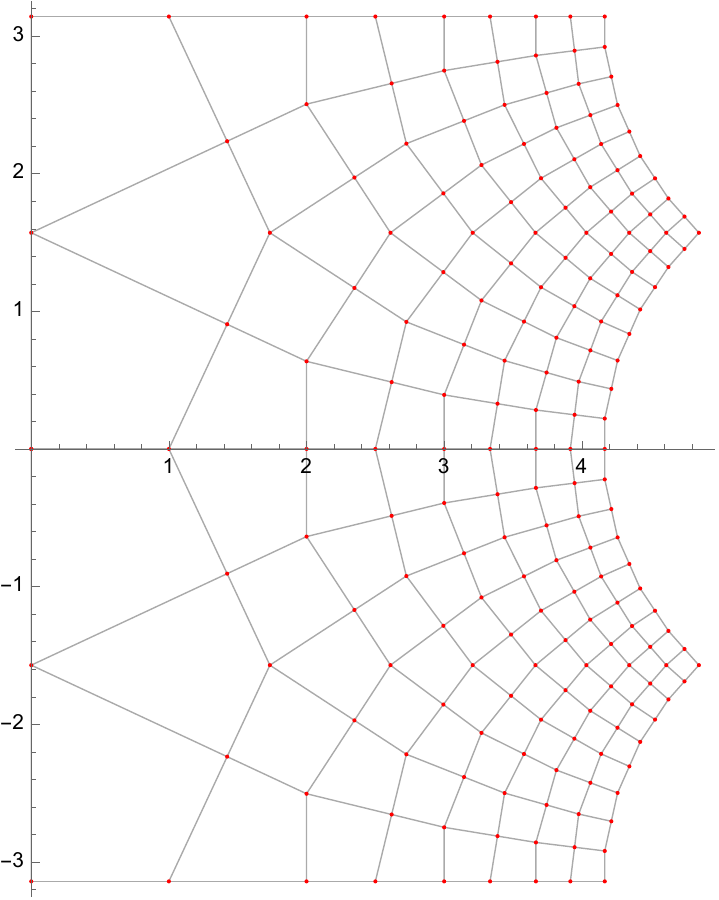}
    \end{minipage}
    \caption{The discrete power function $Z^\frac{2}{5}$ (left) and the logarithmic function $\Log$ (right) with extended domain.}
    \label{fig: powlog_ext}
\end{figure}

The domain of $Z^\gamma$ or $\Log$ can be extended to ``discrete Riemann surface'' by setting initial conditions as follows\cite{ando2014explicit}.
We introduce the polar form of $(m, n)$ as $m+ni=r_{m, n}e^{i\theta_{m, n}}$, and write $f(m, n)=f(r_{m, n}, \theta_{m, n})$ by abuse of notation.
Then, the initial conditions are given by
\begin{gather}
    Z^\gamma (1, k\pi/2)=e^{k\gamma\pi i/2},\\
    \begin{gathered}
        \Log (1, k\pi/2)=k\pi i, \quad \Log (2, k\pi/2)=1+k\pi i,\\
        \Log (\sqrt{2}, \pi/4+k\pi/2)=\left(\frac{\pi}{2}+k\right)i,
    \end{gathered}
\end{gather}
for $k \in \Z$ (see Fig. \ref{fig: powlog_ext}).

\clearpage

\section{Force equilibrium and geometric duality}

\subsection{Basic theory}
\begin{figure}[htbp]
    \centering
    \includegraphics[width=0.5\linewidth]{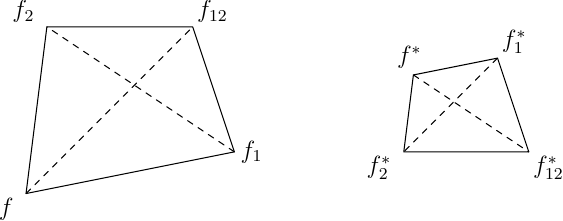}
    \caption{Elementary quadrilateral of $f$ (left) and its dual (right).}
    \label{fig: d_iso-dual}
\end{figure}

A discrete holomorphic function $f$ has a dual net called the \emph{Christoffel dual}.
This duality is the key of the description of equilibrium of truss structures.
\begin{dfn}[Bobenko--Suris\cite{bobenko2008discrete}]
    The \emph{Christoffel dual} of a discrete holomorphic function $f$ is a map $f^* : \Z^2 \to \C$ satisfying
    \begin{equation}
        \begin{gathered}
            f^*_1-f^* \parallel f_1-f, \quad f^*_2-f^* \parallel f_2-f,\\
            f^*_{12}-f^* \parallel f_1-f_2, \quad f^*_1-f^*_2 \parallel f_{12}-f
        \end{gathered}
    \end{equation}
    (see Fig. \ref{fig: d_iso-dual}).
    $f^*$ is unique up to scaling and translation.
\end{dfn}

The Christoffel dual of the discrete power function $Z^\gamma$ is the complex conjugate of $Z^{2-\gamma}$ $(0<\gamma<2)$\cite{bobenko1999discrete}.

To relate these geometric concepts to structural mechanics, let us consider a quadrilateral mesh where internal forces act along the edges at the midpoints.

\begin{thm}[Schief\cite{schief2014integrable}]
    A circular net $\bm{r}$ may be regarded as a discrete membrane in equilibrium with purely tangential internal forces acting at the midpoints of the edges if and only if $\bm{r}$ constitutes a discrete holomorphic function.
    The internal forces are encoded in the Christoffel dual $\bm{r}^*$.
    \label{thm: schief}
\end{thm}
\begin{figure}[htbp]
    \centering
    \includegraphics[width=0.5\linewidth]{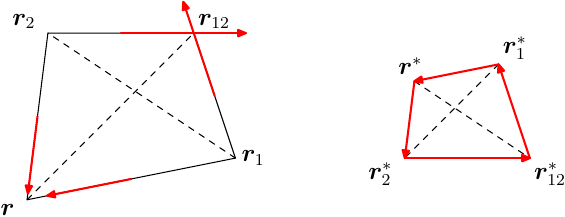}
    \caption{Form diagram (left) and force diagram (right).}
\end{figure}

Here, $\bm{r}$ and $\bm{r}^*$ represent the structure's form and force equilibrium, corresponding to the form and force diagrams in graphic statics, respectively\cite{maxwell1864reciprocal}.

In the following two subsections, we show the method of constructing truss structures in equilibrium\cite{hayashi2024parametric, jikumaru2025pure}.

\subsection{Compression-tension mixed case}

\begin{figure}[htbp]
    \centering
    \includegraphics[width=0.6\linewidth]{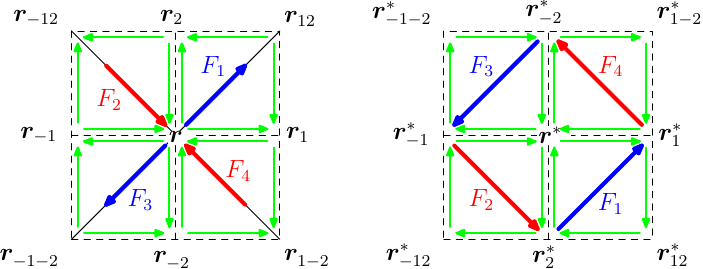}
    \caption{Force equilibrium at node $\bm{r}$.}
    \label{fig: form_force}
\end{figure}

Let $\bm{r}$ be a discrete holomorphic function, and consider four adjacent quadrilaterals meeting at an interior node $\bm{r}$ (see Fig. \ref{fig: form_force}).
The diagonals incident on $\bm{r}$ constitute the truss members.
The axial forces (represented by red or blue arrows) can be decomposed along the edges of $\bm{r}$ (represented by green arrows).
Theorem \ref{thm: schief} implies that the closure of the dual quadrilateral $(\bm{r}^*, \bm{r}^*_1, \bm{r}^*_{12}, \bm{r}^*_2)$ ensures force equilibrium along the member $(\bm{r}, \bm{r}_{12})$.
Furthermore, nodal equilibrium at $\bm{r}$ is equivalent to the closure of the dual cycle $(\bm{r}^*_2, \bm{r}^*_1, \bm{r}^*_{-2}, \bm{r}^*_{-1})$. 
The diagonal length of $\bm{r}^*$ represents the magnitude of axial force acting on the corresponding diagonal member of $\bm{r}$.
Consequently, the relative axial force magnitudes can be uniquely determined by the truss geometry.

\subsection{Pure tensile (compressive) case}

\begin{figure}[htbp]
    \centering
    \includegraphics[width=0.3\linewidth]{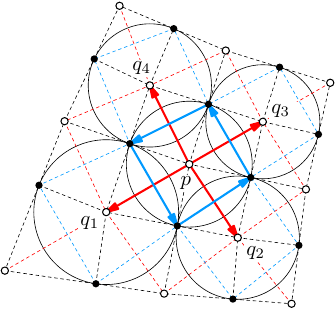}
    \caption{Force equilibrium at node $p$ (pure tension).}
    \label{fig: form_force-pure_tension}
\end{figure}

For the truss structures formed by the diagonals of a discrete holomorphic function under pure tension (or compression), their force diagrams cannot be derived from the diagonals of Christoffel dual, as the forces acting on the corresponding lattice of discrete holomorphic function are not in equilibrium.
However, if the two diagonals of each elementary quadrilateral are orthogonal, the force diagram can be constructed by the following method\cite{jikumaru2025pure}.
As described in \cite{agafonov2000discrete}, all elementary quadrilaterals of the discrete power and logarithmic functions are of kite form, thereby satisfying this orthogonality requirement.
In Fig. \ref{fig: form_force-pure_tension}, the red diagonals of the discrete power (or logarithmic) function (black dotted lines) constitute the truss geometry.
By a $\pi/2$ rotation and scaling of axial force vectors, the blue closed polygon representing the force equilibrium at $p$ is obtained.

\section{Michell--Prager type truss structures}

A \emph{Michell truss} is a truss structure that minimizes its material volume under given boundary conditions\cite{lewinski2018michell}.
Michell effectively reduced this complex, material-dependent optimization problem to a purely geometric problem of truss member arrangement.
\begin{thm}[Michell\cite{michell1904limits}]
    Under given boundary conditions, minimizing the material volume of a truss structure is equivalent to minimizing the \emph{Michell functional}
    \begin{equation}
        \Phi=\sum_i |F_i| L_i,
    \end{equation}
    where $F_i$ and $L_i$ denote the internal (axial) force and the length of the $i$-th member, respectively.
\end{thm}

As a class of truss structures that potentially minimizing $\Phi$, the following property is introduced:
\begin{dfn}[Jikumaru\cite{jikumaru2025pure}]
    A truss structure is called \emph{Michell--Prager type} if there exists a constant $C_0$ such that
    \begin{equation}
        |F_i| L_i = C_0,
    \end{equation}
    for every member $i$.
\end{dfn}
In the compression-tension mixed case described in Section 3.2, truss structures constructed from the discrete holomorphic functions belong to the Michell--Prager type.
For a discrete isothermic net $\bm{r}$, the diagonal lines $(\bm{r}, \bm{r}_{12})$ and $(\bm{r}^*_1, \bm{r}^*_2)$ are related by
\begin{equation}
    \bm{r}^*_1-\bm{r}^*_2=C_0\frac{\bm{r}_{12}-\bm{r}}{|\bm{r}_{12}-\bm{r}|^2}, ~ \bm{r}^*_{12}-\bm{r}^*=C_0\frac{\bm{r}_1-\bm{r}_2}{|\bm{r}_1-\bm{r}_2|^2},
\end{equation}
where $C_0$ is a real constant (see Corollary 4.33 in \cite{bobenko2008discrete}).
This relation implies
\begin{equation}
    |\bm{r}^*_1-\bm{r}^*_2||\bm{r}_{12}-\bm{r}|=|\bm{r}^*_{12}-\bm{r}^*||\bm{r}_1-\bm{r}_2|=|C_0|.
\end{equation}
Since the dual edge lengths $|\bm{r}^*_1-\bm{r}^*_2|$ and $|\bm{r}^*_{12}-\bm{r}^*|$ correspond to the magnitudes of the axial forces, while the diagonal lengths $|\bm{r}_{12}-\bm{r}|$ and $|\bm{r}_1-\bm{r}_2|$ represent the member lengths, the resulting truss satisfies the condition $|F_i|L_i=|C_0|$, thereby confirming its Michell--Prager type property.
We remark that pure tensile (compressive) structures in Section 3.3 cannot be of the Michell--Prager type.

\section{Results of numerical analysis}

\begin{figure}[htbp]
    \begin{minipage}[b]{0.3\linewidth}
        \centering
        \includegraphics[height=3.2cm]{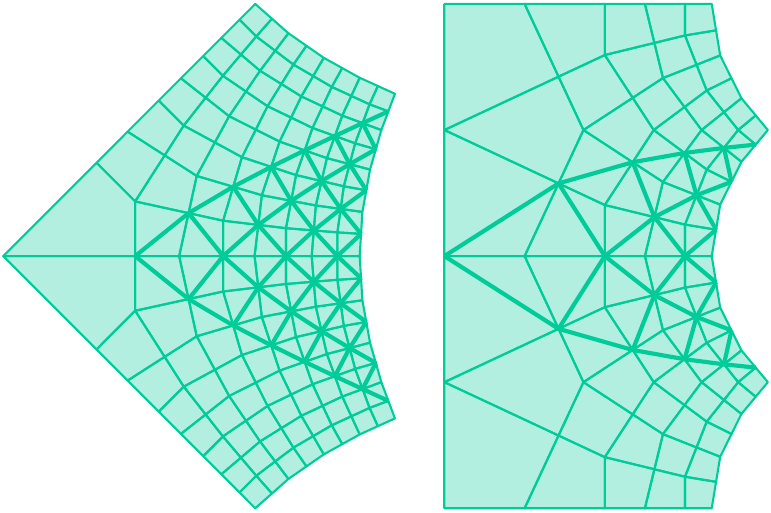}
        \caption{Symmetric type.}
        \label{fig: case-1}
    \end{minipage}
    \hfill
    \begin{minipage}[b]{0.3\linewidth}
        \centering
        \includegraphics[height=3.2cm]{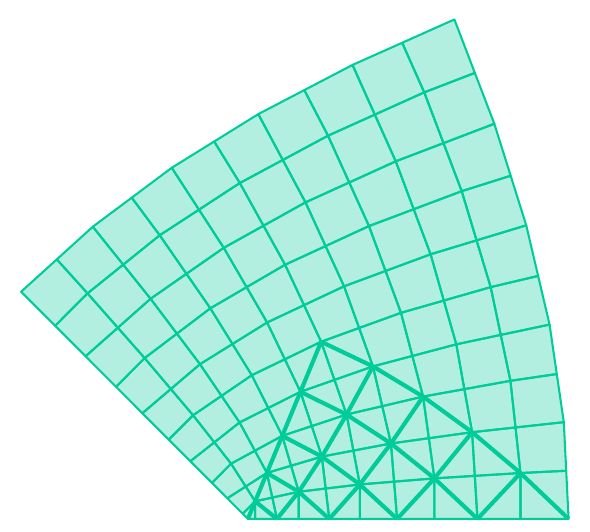}
        \caption{Asymmetric type.}
        \label{fig: case-2}
    \end{minipage}
    \hfill
    \begin{minipage}[b]{0.3\linewidth}
        \centering
        \includegraphics[height=3.2cm]{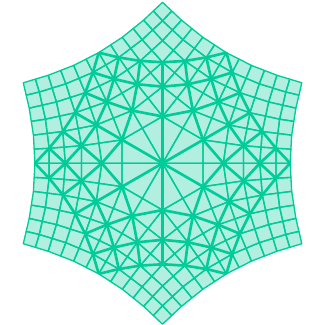}
        \caption{Whole type.}
        \label{fig: case-3}
    \end{minipage}
\end{figure}

In this section, we numerically verify whether the structures derived from the discrete power and logarithmic functions effectively approximate Michell trusses.
Appropriate selection of partial structures and boundary conditions is crucial for generating structures that effectively approximate Michell trusses.
We construct three types of structures: (1) symmetric cantilever type extracted from the image of two quadrants (Fig. \ref{fig: case-1}), (2) asymmetric cantilever type extracted from the image of the first quadrant (Fig. \ref{fig: case-2}), (3) whole type extracted from the image of extended domain (Fig. \ref{fig: case-3}).
In each figure, the bold lines represent the members of the truss structure.
By the optimization of the free node positions minimizing the Michell functional, we measure their deviation from the critical points of the Michell functional.

\noindent
{\bf (1) Symmetric cantilever type structure}

\begin{figure}[htbp]
    \centering
    \begin{minipage}[b]{0.48\linewidth}
        \centering
        \includegraphics[height=3.3cm]{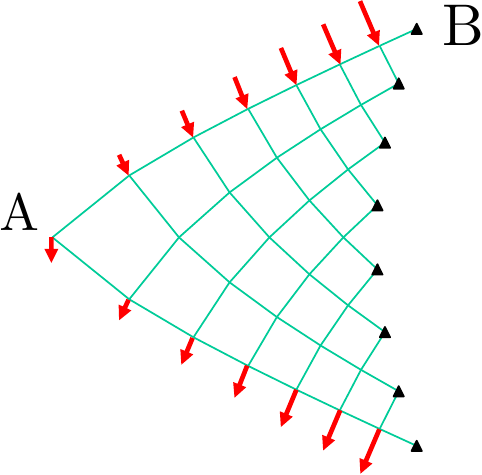}
        \caption{Form diagram.}
        \label{fig: form_1}
    \end{minipage}
    \begin{minipage}[b]{0.48\linewidth}
        \centering
        \includegraphics[height=3.3cm]{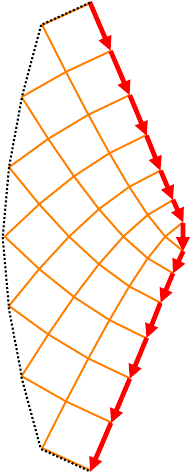}
        \caption{Force diagram.}
        \label{fig: force_1}
    \end{minipage}
\end{figure}

The following figures illustrate the case of $Z^{0.5}$.
The red arrows represent the external loads, and the black triangles mark the fixed supports in Fig. \ref{fig: form_1}.
The relative magnitudes and directions of external loads are determined in order to close the force diagram extracted from the Christoffel dual $Z^{1.5}$ (Fig. \ref{fig: force_1}).

\begin{figure}[htbp]
    \centering
    \begin{minipage}[b]{0.48\linewidth}
        \centering
        \includegraphics[height=3.3cm]{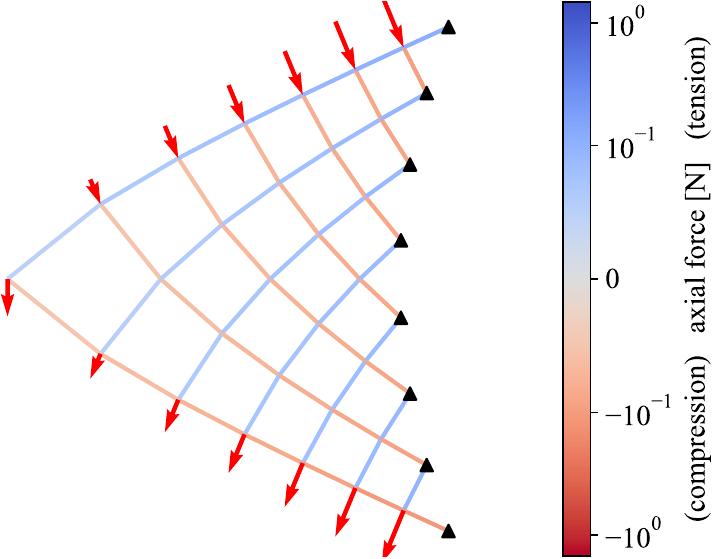}
        \caption{Axial force distribution of the structure.}
        \label{fig: dist_1}
    \end{minipage}
    \hfill
    \begin{minipage}[b]{0.48\linewidth}
        \centering
        \includegraphics[height=3.3cm]{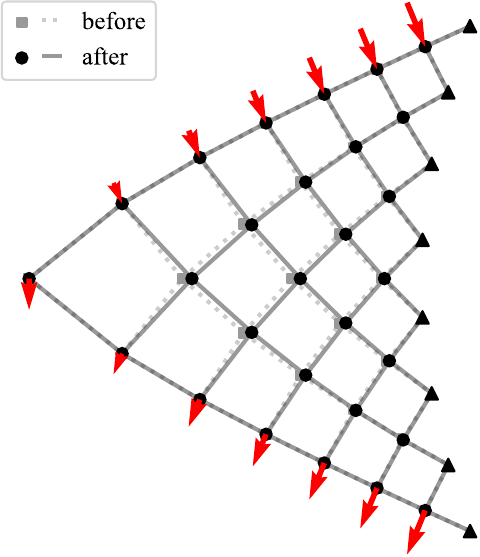}
        \caption{Shape change before and after optimization.}
        \label{fig: opt_1}
    \end{minipage}
\end{figure}

The structure is scaled so as the distance between two points, A and B in Fig. \ref{fig: form_1}, to be 1.
Also, the magnitudes of the loads are normalized so that their sum equals 1.
Figs. \ref{fig: dist_1} and \ref{fig: opt_1} illustrate the results of numerical analysis of the truss structure in Fig. \ref{fig: form_1} under this normalization.
Fig. \ref{fig: dist_1} shows the distribution of compressive (red) and tensile (blue) forces of the members.
The color gradation corresponds to the magnitude of the axial forces.
In Fig. \ref{fig: opt_1}, light gray nodes and black nodes represent the node position before and after optimization, respectively.
It can be observed that the shape remains almost unchanged before and after optimization, which suggests the original shape is approximately the critical point of the Michell functional.

For the optimization process minimizing the Michell functional, the SLSQP method\cite{kraft1988software,scipy_optimize} is used.
The optimization process stops when the change in the objective function value is less than the tolerance of $1.0 \times 10^{-8}$.
Fig. \ref{fig: opt_1} shows the shape change before and after optimization.
Table \ref{tab: pow_0.5} describes the relationship between the number of supports $N_\mathrm{sup}$ (corresponding to the mesh resolution) and the deviation from the critical points of the Michell functional $\Phi$.
Although the maximum displacement of the nodes $\max(\mathbf{D})$ increases, the norm of the gradient $|\nabla \Phi|$ and the objective function values before and after optimization, $\Phi_\mathrm{before}$ and $\Phi_\mathrm{after}$, decrease as $N_\mathrm{sup}$ increases.
This tendency is also observed for other exponents $0<\gamma<1$ and the $\Log$ function.

\begin{table}[htbp]
    \centering
    \begin{tabular}{|ccccc|}
        \hline
        $N_\mathrm{sup}$ & $|\nabla \Phi|$ & $\max (\mathbf{D})$ & $\Phi_\mathrm{before}$ & $\Phi_\mathrm{after}$ \\
        \hline \hline
        $6$ & $6.217 \times 10^{-2}$ & $1.736 \times 10^{-2}$ & $0.607228$ & $0.606542$ \\
        $8$ & $6.102 \times 10^{-2}$ & $1.890 \times 10^{-2}$ & $0.563301$ & $0.562267$ \\
        $10$ & $5.599 \times 10^{-2}$ & $2.096 \times 10^{-2}$ & $0.535163$ & $0.533879$ \\
        $12$ & $5.069 \times 10^{-2}$ & $2.275 \times 10^{-2}$ & $0.515567$ & $0.514101$ \\
        $14$ & $4.594 \times 10^{-2}$ & $2.330 \times 10^{-2}$ & $0.501081$ & $0.499477$ \\
        $16$ & $4.185 \times 10^{-2}$ & $2.393 \times 10^{-2}$ & $0.489894$ & $0.488180$ \\
        \hline
    \end{tabular}
    \caption{Deviation from the critical points of $\Phi$ ($\gamma=0.5$).}
    \label{tab: pow_0.5}
\end{table}

\noindent
{\bf (2) Asymmetric cantilever type structure}

The following figures illustrate the case of $Z^{1.5}$.
The boundary conditions (supports and loads) are set in a manner similar to the previous case.
Consistent with the previous case, higher resolution ($N_\mathrm{sup}$) decreases both $|\nabla \Phi|$ and $\Phi$ values, while increasing $\max(\mathbf{D})$.

\begin{figure}[htbp]
    \centering
    \begin{minipage}[b]{0.48\linewidth}
        \centering
        \includegraphics[width=0.7\linewidth]{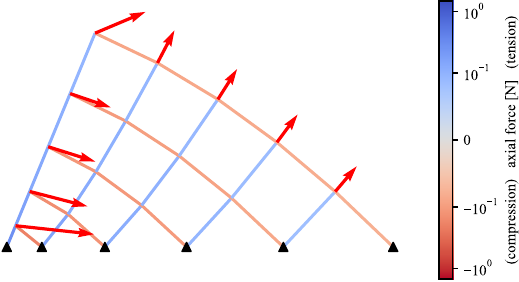}
        \caption{Axial force distribution of the structure.}
    \end{minipage}
    \hfill
    \begin{minipage}[b]{0.48\linewidth}
        \centering
        \includegraphics[width=0.55\linewidth]{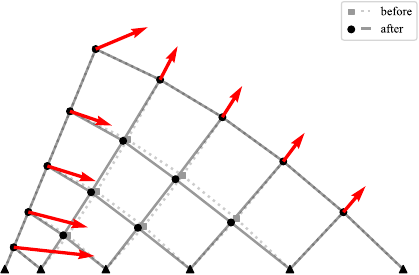}
        \caption{Shape change before and after optimization.}
    \end{minipage}
\end{figure}

\begin{table}[htbp]
    \centering
    \begin{tabular}{|ccccc|}
        \hline
        $N_\mathrm{sup}$ & $|\nabla \Phi|$ & $\max (\mathbf{D})$ & $\Phi_\mathrm{before}$ & $\Phi_\mathrm{after}$ \\
        \hline \hline
        $4$ & $1.025 \times 10^{-1}$ & $1.689 \times 10^{-2}$ & $0.596826$ & $0.595955$ \\
        $5$ & $1.061 \times 10^{-1}$ & $1.764 \times 10^{-2}$ & $0.526496$ & $0.525188$ \\
        $6$ & $1.023 \times 10^{-1}$ & $1.984 \times 10^{-2}$ & $0.485503$ & $0.483961$ \\
        $7$ & $9.751 \times 10^{-2}$ & $2.084 \times 10^{-2}$ & $0.458567$ & $0.456882$ \\
        $8$ & $9.289 \times 10^{-2}$ & $2.101 \times 10^{-2}$ & $0.439446$ & $0.437666$ \\
        $9$ & $8.872 \times 10^{-2}$ & $2.224 \times 10^{-2}$ & $0.425124$ & $0.423275$ \\
        \hline
    \end{tabular}
    \caption{Deviation from the critical points of $\Phi$ ($\gamma=1.5$).}
    \label{tab: pow_1.5}
\end{table}

The results of \textbf{(1)} and \textbf{(2)} suggest that these Michell--Prager type structures approximate Michell trusses well.

\newpage

\noindent
{\bf (3) Whole type structure}

\begin{figure}[htbp]
    \centering
    \begin{minipage}[b]{0.48\linewidth}
        \centering
        \includegraphics[height=4cm]{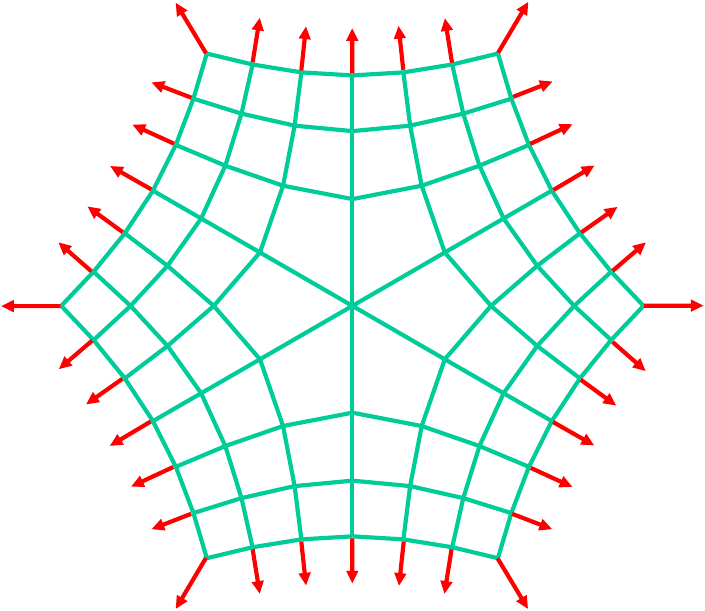}
        \caption{Form diagram.}
        \label{fig: form_23}
    \end{minipage}
    \hfill
    \begin{minipage}[b]{0.48\linewidth}
        \centering
        \includegraphics[height=4cm]{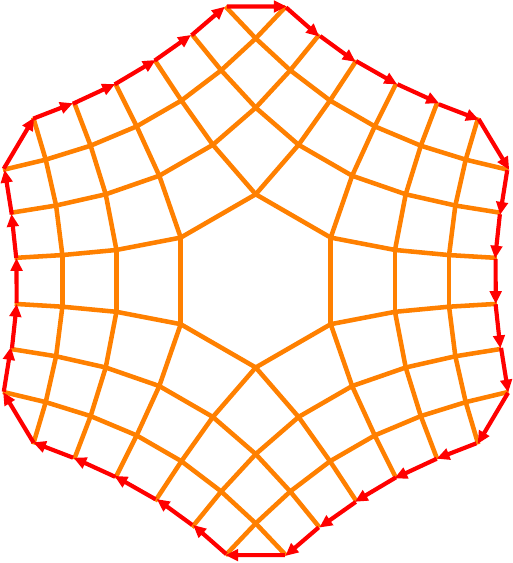}
        \caption{Force diagram.}
        \label{fig: force_23}
    \end{minipage}
\end{figure}

\begin{figure}[htbp]
    \centering
    \begin{minipage}[b]{0.48\linewidth}
        \centering
        \includegraphics[height=4cm]{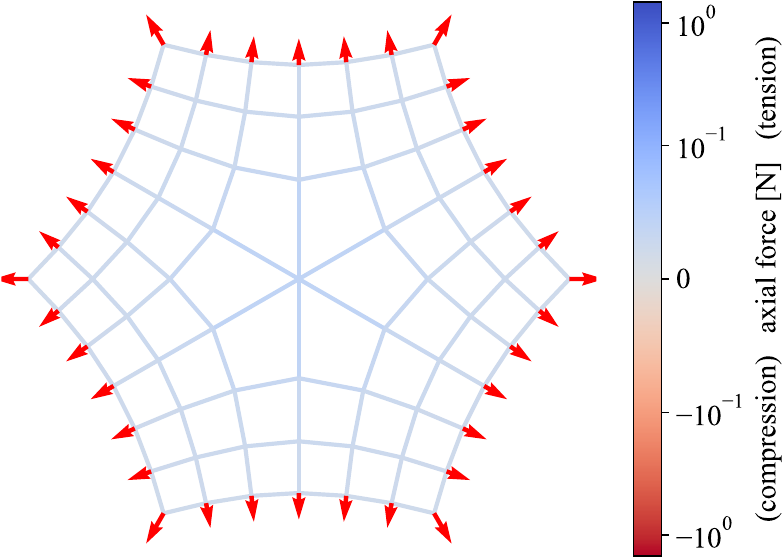}
        \caption{Axial force distribution of the structure.}
        \label{fig: dist_23}
    \end{minipage}
    \hfill
    \begin{minipage}[b]{0.48\linewidth}
        \centering
        \includegraphics[height=4cm]{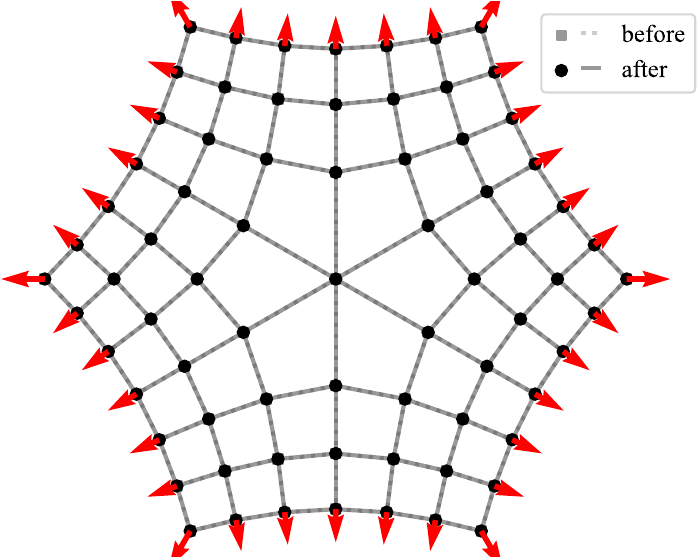}
        \caption{Shape change before and after optimization.}
        \label{fig: opt_23}
    \end{minipage}
\end{figure}

In the case of some exponents $\gamma$ satisfying $4/\gamma \in \N$, whole type pure tension structures can be obtained (see Figs. \ref{fig: form_23} and \ref{fig: force_23} for specific boundary conditions).
Remarkably, these structures remain invariant during the optimization process, despite not being of the Michell--Prager type.
This finding implies these structures numerically coincide with the critical points of the Michell functional $\Phi$, strongly suggesting that they are Michell trusses.
Figs. \ref{fig: dist_23} and \ref{fig: opt_23} illustrate the results for $Z^\frac{2}{3}$.

\section*{Acknowledgments}
This work was supported by JST CREST Grant No. JPMJCR1911 and JSPS KAKENHI No. 25K21661.

\end{document}